\documentclass[16pt]{article}
\usepackage{graphicx}
\usepackage{amssymb}
\usepackage{color}

 \let\cal\mathcal

\usepackage{graphicx}
\usepackage[curve]{xypic}
\usepackage{tikz}

\newcommand\C{{\mathbb C}}
\newcommand\Z{{\mathbb Z}}

 \newtheorem{theorem}{Theorem}[subsection]
\newtheorem{proposition}[theorem]{Proposition}
\newtheorem{corollary}[theorem]{Corollary}
\newtheorem{lemma}[theorem]{Lemma}
\newtheorem{definition}[theorem]{Definition}
\newtheorem{example}[theorem]{Example}
\newtheorem{remark}[theorem]{Remark}

\usepackage{graphicx}
\begin{document}
\title{ The Topology of the Normalization of Complex  Surface Germs  }
\author{Francoise Michel}

\maketitle

\begin{abstract}

 Let  $(X,p)$ be a  reduced   complex surface germ and  let $L_X$ be its well defined link. If $(X,p)$ is normal at $p$,  D. Mumford \cite{mu} shows that $(X,p)$ is smooth if and only if $L_X$ is simply connected.  Moreover, if $p$  is an isolated singular point, $L_X$ is a three dimensional  Waldhausen graph manifold. Then, the    Plumbing Calculus  of 
W. Neumann \cite{ne}  shows  that the homeomorphism class of $L_X$ determines a unique   plumbing in  normal form  and consequently,  determines   the topology of the good minimal resolution of $(X,p)$.

Here, we do not assume that  $X$ is normal at $p$, and so,    the singular locus  $(\Sigma ,p)$  of $(X,p)$ can   be one dimensional.  We describe the topology of the singular link $L_X$ and we show that the homeomorphism class of $L_X$  (Theorem \ref{th1}) determines the homeomorphism class of the normalization and consequently the plumbing of  the minimal good resolution of $(X,p)$. 
Moreover, in  Proposition \ref{pr1}, we   obtain  the  following generalization   of the  D. Mumford  theorem \cite{mu}:

 Let $\nu  :  (X',p')  \to  (X,p)$  be the normalization of   an irreducible  germ of complex surface $(X,p)$. If the link $L_X$ of $(X,p)$ is simply connected then  $\nu$ is a homeomorphism and $(X',p')$ is a  smooth germ of surface.

\end{abstract}



 \noindent {\bf Adresses.}

\noindent Fran\c coise Michel / Laboratoire de Math\' ematiques Emile
Picard  /
Universit\' e Paul Sabatier / 118 route de Narbonne / F-31062
Toulouse / FRANCE

e-mail: fmichel@picard.ups-tlse.fr

 \bigskip
 \section{Introduction}

 Let $I$ be a reduced ideal in $\C \{z_1, \dots ,z_n \} $ such that the quotient algebra $A_X=\C \{z_1, \dots ,z_n \} / I$ is two-dimensional.
The zero locus,  at the origin  $0$ of  $\C ^{n}$, of a set of generators of $I$ is an analytic surface germ  embedded in $(\C^{n},0).$   Let $(X,0)$ be its intersection with  the  compact   ball $B_{\epsilon}^{2n}$ 
of radius a sufficiently small  $\epsilon$,  centered at the origin in $\C^{n}$,  and let  $L_X$  be its  intersection with the boundary   $S_{\epsilon} ^{2n-1}$ of $B_{\epsilon}^{2n}.$  Let $\Sigma $ be the set of the singular points of $(X,0).$ 

As $I$ is reduced,  $\Sigma$ is empty when $(X,0)$ is smooth, $\Sigma$  is equal to the origin when $0$ is an isolated singular point and  $\Sigma$  is a curve when the germ has a non-isolated singular locus (in particular $(X,0)$  can be a   reducible germ). 

 If $\Sigma $ is a curve,   $K_{\Sigma }=\Sigma \cap S_{\epsilon} ^{2n-1}$ is the  disjoint union of $r$ one-dimensional circles ($r$ being the number of irreducible components of $\Sigma$) embedded in $L_X$. We  say that $K_{\Sigma }$ is the link  of $\Sigma.$ By the conic structure theorem of J.  Milnor \cite{mi68}, for a sufficiently small $\epsilon$,  $(X,\Sigma, 0)$ is homeomorphic to the cone on  the pair $(L_X, K_{\Sigma }) $.   When $\Sigma =\{ 0 \} $,  $(X, 0)$ is homeomorphic to the cone on $L_X$.

 On the other hand,  thanks to A.Durfee   \cite{du},  the  homeomorphism class of $(X,0)$ depends only on the isomorphism class of the algebra $A_X$ (i.e. is independent of  a sufficiently small $\epsilon$ and of the choice of the embedding in $(\C ^n,0)$).  
The analytic type of $(X,0)$ is given by the isomorphism class of $A_X$,  and  its  topological type  is given by the homeomorphism  class of  $(X, 0)$. 

\begin{definition}  {\bf The  link of $(X,0)$}  is the homeomorphism class of   $L_X$.  {\bf The  link of $(\Sigma,0)$} is the homeomorphism class of the pair $(L_X, K_{\Sigma }).$

 \end{definition}
 
  Let $\nu  :  (X',p') \rightarrow (X,p)$  be the normalization of   a reduced germ of  complex surface $(X,p)$.

\begin{remark}\label{remnor}  
 \begin{enumerate}

 \item If $(X,0)$ is reducible, let $(\cup _{1\leq i \leq r}X_i,0)$ be its decomposition as a union of irreducible surface  germs. Let $\nu_i : (X'_i,p_i) \to (X_i,0)$ be the normalization of the irreducible components of $(X,0)$. The  morphisms $\nu_i$ induce the   normalization morphism   on the disjoint union $\coprod _{1\leq i \leq r}(X'_i,p_i)$. 
 
 \item If $0$ is an isolated singular point of $(X,0)$ then $\nu $ is a homeomorphism. 
 
 \end{enumerate}

\end{remark} 

\begin{definition}\label{de1}

 If $(\Sigma ,0)$ is a one-dimensional germ, let $\sigma$ be an irreducible component of $\Sigma $. Let $\sigma '_j , 1\leq j \leq n(\sigma), $ be the  $n(\sigma)$  irreducible components of $\nu ^{-1} (\sigma)$   and   let $d_j$ be the degree of $\nu $ restricted to $\sigma'_j $.  The following number $k(\sigma)$:
 $$ k(\sigma)=:d_1+\dots +d_j+\dots +d_{n(\sigma )}.$$
 is {\bf{the total degree of $\nu$ above $\sigma$}}.
\end{definition}

 Let $\Sigma_+$ be the union of the irreducible components $\sigma $ of  $\Sigma $ such that $k(\sigma )>1$. In $L_X$,  let $K_{\Sigma_+}$ be the link of $\Sigma_+$.  We choose  a compact regular neighbourhood $N(K_{\Sigma_+})$ of $K_{\Sigma_+}$. Let $E(K_{\Sigma_+})$ be the closure of $L_X \setminus N(K_{\Sigma_+})$. By definition  {\bf {$E(K_{\Sigma_+})$ is the (compact) exterior of  $K_{\Sigma_+}.$ }}

  In  Section  3 of \cite{mi},  one can find a description of  the topology of  $N(K_{\Sigma_+})$ which implies the following lemma \ref{lem1}. To be  be self-contained, we begin   Section 3  with a quick proof of it.
 
 {\bf{Lemma}  (\ref{lem1})} \\
   1.  The restriction of $\nu $ to $\nu ^{-1} (E(K_{\Sigma_+})) $ is an homeomorphism and  $(L_X\setminus K_{\Sigma_+})$ is a topological manifold.\\
 2.  The link   $K_{\Sigma_+}$ is the set of the topologically singular points of $L_X$.\\
   3.   The homeomorphism class of $L_X$  determines  the homeomorphism class of
   $N(K_{\Sigma_+})$ and  $E(K_{\Sigma_+})$. \\
    4.  The number of connected components of $E(K_{\Sigma_+})$ is equal to the number of irreducible components of $(X,0).$

In Section 3, we prove the following theorem:

{\bf {Theorem} (\ref{th1})}
 Let $(X,0)$ be a reduced surface germ. The homeomorphism class of $L_X$  determines the homeomorphism class of the link $L_{X'}$ of the normalization of $(X,0)$.

 The proof of   Theorem \ref{th1} gives an explicit construction to transform $L_X$ in a Waldhausen graph manifold orientation preserving homeomorphic to $L_{X'}$.  Then,    the plumbing calculus of W. Neumann \cite{ne}  implies the following corollary.
 
 \begin{corollary}\label{co1} 
 The homeomorphism class of $L_X$  determines the topology of the plumbing of the good minimal resolution of  the germ $(X,0)$.
 
 \end{corollary}
 
 The proof of \ref{th1} is based on a detailed description of a regular neighbourhood $N(K_{\Sigma_+})$ of the topologically  singular locus $K_{\Sigma_+}$  of  $L_X$ and on the topology of $\nu$ restricted to $\nu ^{-1} ( L_X \setminus K_{\Sigma_+}) $. 
 
  In Section 2, we  describe  the topology of  a {\bf{$d$-curling }} and the topology of a  {\bf{ singular pinched torus}} which is defined as   the mapping torus  of  an orientation preserving homeomorphism acting on a reducible germ of curves. Curlings and singular pinched tori are already  studied in \cite{mi}. But, to prove  the new results \ref{th1}, \ref{lem2} and \ref{pr1} of this paper, we need to insist on particular properties, given  in  \ref{remcurl},  of these two topological objects.
      Section 2 contains  also  the presentation of  the following example which is a typical  example of   $d$-curling.

{\bf{Example}} (\ref{ex1})
Let $X= \{(x,y,z)\in \C ^3 \ where \  z^d-xy^d=0  \ and\ d>1 \} $. The normalization of $(X,0)$ is smooth i.e.  the morphism   $\nu : (\C^2,0) \to (X,0)$  defined  by $(u,v) \mapsto  (u^d,v,uv)$ is a normalization morphism.  The singular locus of $(X,0)$ is the line $ l_x=\{ (x,0,0) \in \C^3, x\in \C \}$. Let $T= \{(u,v)\in  (S \times D ) \subset \C^2 \}$.      We have  $n(l_x)=1$ because $\nu^{-1}(l_x ) $ is the line $ \{(u,0) \in \C^2, u\in \C \}$ and $d_1=d$. Moreover  $N(K_{l_x}) = \nu (T)$ is   a tubular neighbourhood of $K_{l_x}$ and $\nu $ restricted to $T$ is  the quotient called $d$-curling. In this example $L_X$ is not simply connected. In fact:  $H _1 (L_X,\Z)=\Z / d\Z$.

  In  Section 3, we show that  each connected component  of $N(K_{\Sigma_+})$  is a singular pinched torus (see also  \cite{mi}).  As stated in \cite{l-p} and as proved in \cite{mi}, It implies that   $\nu$ restricted to $\nu ^{-1} ( L_X \setminus K_{\Sigma_+}) $    is   the composition of two kind of quotients: curlings and identifications.
 Here, we need the  more  detailed  description of $N(K_{\Sigma_+})$ given Section 2,  to prove the following Lemma \ref{lem2} which shows how the topology of $L_X$ determines the invariants    $ n(\sigma )$ and $  d_j, 1\leq j \leq n(\sigma ),$  of the normalization morphism.
More precisely:  

{\bf {Lemma (\ref{lem2}) }}  Let  $\sigma$ be an irreducible component of $\Sigma_+ $,   let $K_{\sigma}$ be the link of $\sigma $ in $L_X$. We choose,  in $L_X$,  a compact regular neighbourhood $N(K_{\sigma})$ of $K_{\sigma}$.  
The link    $K_{\sigma}$ is a deformation retract of  $N(K_{\sigma})$.  If  $l_{\sigma } $ is  the homotopy class of  $K_{\sigma}$  in $\pi _{1} (N(K_{\sigma }) )$,  
  then $\pi _1 (N(K_{\sigma}) = \Z  .  l_{\sigma}$. We have:\\
 1.  The tubular neighbourhood   $ \nu ^{-1}(N(K_{\sigma}))$ of $ \nu ^{-1}(K_{\sigma})$ is the disjoint union of $n(\sigma )$ solid tori \\ $T'_j,  1\leq j \leq n(\sigma )$,  and 
the boundary of  $N(K_{\sigma})$ is the disjoint union of $n(\sigma )$  tori.  \\
 2. Let $c$ be the permutation of $k(\sigma )$ elements which is the composition  of $n(\sigma)$ disjoint cycles $c_j$ of order $d_j$.  Then  $N(K_{\sigma})$ is homeomorphic to a singular pinched torus $T(k(\sigma)(D),c)$ wich   has $n(\sigma)$ sheets $T_j$ where  
  $T_j=\nu(T'_j) , 1\leq j \leq n(\sigma)$.\\
3.  On  each connected component $\tau_j$ of the boundary of $N(K_{\sigma})$, the homeomorphism class of $N(K_{\sigma})$ determines a  unique (up to isotopy)  meridian curves $m_j$.  If $l_j$ is a  parallel on  $\tau_j$ the homotopy class of $l_j$ in $\pi _{1} (N(K_{\sigma }) )$ is equal to $d_j.l_{\sigma }$.

  In \cite{mu}, D. Mumford proves that a normal surface germ which has a simply connected link is  a smooth germ of surface. However,  there exist  surface germs   with one dimensional singular locus   and  simply connected links. Obvious examples are  obtained as follows: 
  
    Let $ f(x,y)$ be an irreducible element of $ \C \{x,y\}$ of multiplicity   $m>1$ at $0$.\\  Let $Z= \{(x,y,z)\in \C ^3 \ such \ that  \  f(x,y)=0\}$.  The singular locus of $(Z,0)$ is   the line  $ l_z= \{(0,0,z), \ z\in \C \}$. The  normalization   $ \nu: (\C^2,0) \to (Z,0)$ can be given by a Puiseux expansion of $f(x,y)$. So, the  link $L_Z$ is the sphere  $S^3$.   L\^e 's conjecture states   that this family  it is the only family  of  singular irreducible  surface germs  with one dimensional singular locus and simply connected links (see \cite{l-p} and \cite{bo} for  partial results).

   In Section 4, we prove the following proposition which is a kind of generalization of Mumford's theorem   for non-normal surface germs:
 
 {\bf{Proposition (\ref{pr1} )}} Let $(X,0)$ be an irreducible surface germ. If the link $L_X$ of $(X,0)$ is simply connected then the normalization $\nu : (X',p') \to (X,0)$ is a homeomorphism and $(X',p')$ is smooth at $p'$. In particular,  the normalization is the good minimal resolution of  $(X,0)$.

 \begin{remark}\label{rem1}    There exist reducible  surface germs with simply connected link  for which  the normalization is not a homeomorphism. 
  For example let $Y= \{(x,y,z)\in \C ^3 \ where \  xy=0\}$. \\ 
    In $\C^3$, $(Y,0)$ is the union of two planes and $\Sigma =\Sigma_+= l_z= \{(0,0,z), \ z\in \C \}$.   Moreover,    the normalization is the obvious quotient  $\nu: (\C^2,0 )\coprod (\C^2 ,0) \to (Y,0)$ given by $\nu (x,z)=(x,o,z)$ and  $\nu (y,z)=(0,y,z)$.  Let $ K_i \subset  S^3_i, i=1,2,$  be  two copies of $S^3$ given with  a trivial knot  $K_i$. The link $L_Y$ of  $(Y,0)$ is the quotient  of the disjoint union  $S^3_1 \coprod S^3_2$,    by the identification point by point of  the   trivial knots $K_1$ and $K_2$.  So $L_Y$ is simply connected. But   $(L_Y \setminus N(K_{l_z}))$ is not connected because it  is the disjoint union of two solid tori. Moreover,  a regular neighbourhood  $N(K_{l_z})$,  of the link of the singular locus $l_z$ of $(Y,0)$,    is a pinched singular torus $T(2(D), id)$ as defined in \ref{deperm}, i.e. $N(K_{l_z})$ is the mapping torus  of the identity acting on  a $2$-pinched disc.
 
 \end{remark}

 \subsection{ Conventions}

The boundary of a topological manifold $W$ will be denoted by $b(W).$

A {\bf disc}  (resp. an {\bf  open disc}) will always be an oriented   topological  manifold  orientation preserving homeomorphic  to $ D=\{ z\in \C ,\vert z \vert \leq 1\}$ (resp. to $ \dot D= \{ z\in \C ,\vert z \vert < 1\}$). 
A {\bf circle } will always be an oriented     topological  manifold orientation preserving homeomorphic  to $ S  = \{ z\in \C ,\vert z \vert =1\}.$

\bigskip
{\bf Acknowledgments: } 
 I thank Claude Weber for  useful discussions and  for making many  comments about the redaction of this text.
 
 \vskip.1in

\section{The topology of   $d$-curlings }

 In this section we study in details the topological properties of   $d$-curlings  because  they are  the key to make the proofs of the two original statements of this paper self-contained. In particular, we need well defined, up to isotopy,   meridian curves on the boundary of curlings  since  the proof of   Theorem  \ref{th1}  is based on Dehn fillings associated to  these  meridian curves.  In this section ,  we suppose that $d>1$ to avoid the trivial case $d=1$.

\begin{definition}\label{decurl}  

\begin{enumerate}

 \item  A  {\bf $d$-curling  ${\cal {C}}_d$} is  a topological space  homeomorphic to    the following quotient of a solid torus $S\times D:$
 $$ {\cal {C}}_d=S\times D /(u, 0)\sim (u',0) \Leftrightarrow u^d=(u')^d.$$
 Let $q:(S\times D)\to {\cal {C}}_d$  be the associated quotient morphism. By definition,   $l_0=q(S\times \{0\})$ is the {\bf core of ${\cal {C}}_d$}.  By definition $ {\cal {C}}_d$ is given with the following orientations: The oriented circle $S$ and the oriented disc $D$ induce an orientation on the circles $l_0$, $q(S\times \{z\}),\ z\in D$,  and on the topogical discs $q(\{u\} \times D),\ u\in S$.
 
 \item A  {\bf {$d$-pinched disc}},   $d(D)$,  is   orientation preserving homeomorphic to  the quotient of the disjoint  union of  $d$  oriented  and ordered discs $D_i, 1\leq  i\leq d$ with  origin $0_i$,  by the relation $0_i \sim 0_j$  for all  $i,   1\leq i \leq d$,  and $j, 1\leq j \leq d$. So, all the origins $0_i,  \ 1\leq i  \leq d \ $ are identified in a unique point $\tilde0 $.  By definition $\tilde0 $ is the {\bf {origin of $d(D)$}}. The class  $\tilde{D}_i$ of each disc $D_i$ in $d(D)$ is {\bf an irreducible component } of $d(D)$.
 
 \end{enumerate}
\end{definition}
   
 The kernel of  the homomorphism  $i_1 : \pi_1 (S\times S)\to    \pi_1 (S\times D))$ induced by the inclusion $(S\times S)\subset (S\times D)$ is infinite cyclic generated by the class of the closed simple curve 
  $$m=(\{u\} \times S) =b(\{u\} \times D), u\in S.$$ 
   Moreover, $m$ is oriented as the boundary of the  oriented disc $(\{u\} \times D)$.  This defines a  unique generator $m^+$ of the kernel of $i_1$. So,   any  closed simple curve on the boundary  of $(S\times D)$  which generates the kernel of $i_1$,  can be oriented to be    isotopic to $m$. By definition it is a  meridian curve of $(S\times D)$.  The $d$-curling  $ {\cal {C}}_d$ is defined by the quotient   $q:(S\times D)\to {\cal {C}}_d$. But, $q$ restricted to $(S\times S)$ is the identity. Moreover  $q (\{u\} \times D)$ is a topological disc in $ {\cal {C}}_d$.
 
  So, the kernel of  the homomorphism  $\bar{i}_1 : \pi_1 (b( {\cal {C}}_d))\to    \pi_1 ( {\cal {C}}_d))$ induced by the inclusion $b( {\cal {C}}_d) \subset ( {\cal {C}}_d)$ is infinite cyclic generated by   the class $m^+$ of the  oriented  simple closed  curve  $q(m)$. As for the solid torus, any  simple closed curve in  $b( {\cal {C}}_d)$  which generates the kernel of $\bar{i}_1 $ can be oriented to be isotopic to $q(m)$.

\begin{definition}\label{demer}  An oriented     simple closed curve  $m$,  on the boundary of a $d$-curling  ${\cal {C}}_d$,  is  a {\bf {meridian curve of ${\cal {C}}_d$}} if   the class of $m$,   in  the kernel of $\bar{i}_1$,  is equal to $m^+$.   Let $m$ be a meridian curve of ${\cal {C}}_d$. An   oriented  simple closed curve  $l $ on the boundary of ${\cal {C}}_d$ is a  {\bf {parallel  curve of ${\cal {C}}_d$}} if  $m \cap l=+1$.

\end{definition}

We gather together  the topological properties of a  $d$-curling  ${\cal {C}}_d$ in the following remark.

\begin{remark}\label{remcurl}
\begin{enumerate}

\item  A consequence of the definition \ref{demer} is:
 A meridian curve on  the torus  $b( {\cal {C}}_d)$ is unique up to isotopy and depends only on the homeomorphism class of $ {\cal {C}}_d$.\\ A simple closed curve $\gamma $ on  a torus  $\tau$  is  essential if $\gamma $ does not bound a disc in $\tau$. But   $b( {\cal {C}}_d)$ is a torus and, by  definition,  a meridian curve $m$ of  $ {\cal {C}}_d$ is a generator of $\pi_1 (b( {\cal {C}}_d))$. So, $m$ is essential on $b( {\cal {C}}_d)$.
 \\ Let $(u,z)\in S\times S$,  we have seen that $m=q (\{u\} \times D)$ is a meridian curve of ${\cal {C}}_d$.  If $z\in S=b(D)$, let  us consider $l_z=q( S\times \{z\} )$. So,  we  have $m\cap l_z=+1$ and $l_z$ is a parallel curve of ${\cal {C}}_d$. But,  contrary to meridians, parallels are not unique up to isotopy.
 
 \item A $d$-curling  ${\cal {C}}_d$ can be retracted by deformation onto  its core $l_0$.  Let $l^+$ be the class of $l_0$ in $\pi _1({\cal {C}}_d)$.  If $z\in  (D \setminus \{0 \} ) $, the homotopy class of $q(S\times \{z\}) $ is equal to $d. l^+.$ Moreover, whatever the parallel curve we choose, its class in $\pi _1({\cal {C}}_d)$ is always equal to $d.l^+$.

\item A germ of complex curve $(\Gamma ,p)$ with $d$ irreducible components is a $d$-pinched disc  and $p$ is its  origin.
\item Let $q:(S\times D)\to {\cal {C}}_d$ be a $d$-curling. Let $\pi _d : (S\times D) \to  (S\times D) $ be the covering of degree $d$ defined,  for  all  $ (u,z)\in (S\times D)$,  by $\pi _d (u,z)=(u^d,z)$.   But, $\pi _d$ induces a unique  topological morphism  $\bar{ \pi} _d :  {\cal {C}}_d \to  (S\times D) $ such that  $\pi _d=\bar{\pi} _d \circ q$. By construction, for all $t\in S$,   $ {\cal{D}}_t=\bar{\pi}_d^{-1} (\{t\} \times D) $ is a $d$-pinched disc. If $u^d=t$, the origin of $ {\cal{D}}_t$  is    $q(u,0)$. 

\item   The circles $(S \times \{z\}), \ z\in D,$ equip the solid torus $T=(S\times D)$ with a  trivial fibration  in oriented circles. If we choose   $t\in S$, the first return  map along these circles induces the identity on the disc $(\{t\}\times D)$. Using $\bar{\pi}_d^{-1}$,  we can lift these  fibration by circles on $ {\cal{C}}_d$. Let $h$ be  the   automorphism of $ {\cal{D}}_t$  defined by the first return map  along these circles. So, $h$  is an orientation preserving homeomorphism of $ {\cal{D}}_t$ which induces    a cyclic  permutation  of the $d$ irreducible  components of $ {\cal{D}}_t $. Obviously $h$  keeps the origin of  $ {\cal{D}}_t $ fixed. So, the  $d$-curling  ${\cal {C}}_d$ is the mapping torus of  an orientation preserving homeomorphism which induces  a cyclic permutations of the  $d$ irreducible components of $ {\cal{D}}_t$.

\item As $d>1$, a  homeomorphism of a $d$-pinched disc keeps always the origin fixed. There is, up to isotopy, a unique  orientation preserving  homeomorphism which induces a cyclic permutation of the $d$ discs  irreducible components of  a $d$-pinched disc. So, the mapping torus of   an orientation   preserving homeomorphism  which induces a cyclic permutation of the $d$ irreducible components of a $d$-pinched disc,   is always orientation preserving homeomorphic to a $d$-curling.
\end{enumerate}
\end{remark}

The following example illustrates Definitions \ref{decurl} and Remarks \ref{remcurl}.

\begin{example}\label{ex1}
Let $X= \{(x,y,z)\in \C ^3 \ where \  z^d-xy^d=0\}$. The normalization of $(X,0)$ is smooth i.e. $\nu : (\C^2,0) \to (X,0)$ is given by $(u,v) \mapsto  (u^d,v,uv)$.  Here $B= \nu (D\times D)$  is a good   semi-analytic neighbourhood of $(X,0)$ in the sense of A. Durfee \cite{du}.  So, the link of $(X,0)$ can be  defined as  $L_X= X\cap \nu ( (S\times D) \cup (D\times S) )$. Let $T= \{(u,v)\in  (S \times D ) \subset \C^2 \}$.   Let $\pi_x : \nu (T) \to   S $ be the projection  
$(x,y,z) \mapsto   x$ restricted to $\nu (T)$. Here the singular locus of $(X,0)$ is the line $ l_x=\{ (x,0,0) \in \C^3, x\in \C \}$ and   $N(K_{l_x}) = L_X \cap (\pi_x^{-1}(S)) =\nu (T)$ is  a tubular neighbourhood of $K_{l_x}$. 

 Let   $q: T \to \cal {C}_d$ be the quotient morphism  which defined a $d$-curling (see \ref{decurl}).
There exists  a well defined homeomorphism  $f: \cal {C}_d \to N(K_{l_x})$ which satisfies 
$f(q(u,v))=(u^d,v,uv)$. So,  $N(K_{l_x})$ is a    d-curling  and  $K_{l_x}$ is its core. Moreover, $f$ restricted to the  core $l_0$  of $ \cal {C}_d $  is a homeomorphism onto  $K_{l_x}$.

Let us take  $s=e^{2i\pi /d}$.  The intersection $ {\cal{D}}_1=N( K_{l_x} )\cap \{x=1\}=  \{(1,y,z)\in \C ^3 \ where \  z^d-y^d=0\}$ is a plane curve germ at $(1,0,0)$ with  $d$ irreducible components  given by  $  \nu (s^k\times D),\ 1\leq k \leq d$. \\
 On the torus $\tau =b(N( K_{l_x} ))= \nu (S\times S)$,   $m=\nu (\{1\} \times S)$ is a meridian curve of $N( K_{l_x} )$  and $l_1= \nu(S\times \{1\})$ is a parallel.  Moreover, $N( K_{l_x} )$ is saturated by the foliation in  oriented circles $l_v=  \nu (S\times \{v\})$ which cuts  $ {\cal{D}}_1$ transversally at  the $d$ points $\{ (1,v, s^kv), \ 1\leq k \leq d \}, $ when  $v\neq 0$ and at $(1,0,0)$ when $v=0.$ So, $N( K_{l_x} )$ is the mapping torus of  the homeomorphism defined  on  the  $d$-pinched disc $ {\cal{D}}_1$   by the first return map   along the circles $l_v$.

To compute $H_1(L_X,\Z )$, we use the Mayer-Vietoris  sequence associated to the decomposition of $ L_X$ as the union $ N( K_{l_x} ) \cup \nu (D\times S) $.  The homology  classes 
 $\bar{m}$ and $\bar{l}_1$ of the  curves $m$ and $l_1$ form  a basis of    $H_1(b(N( K_{l_x} )) ,\Z)$. But, $\bar{m}$
  is a generator of $H_1(\nu (D\times S),\Z)$ and is equal to $0$ in $H_1(N( K_{l_x} ) ,\Z)$. 
  As the class  $\bar{l}_0$ of $K_{l_x}$ is a generator of $H_1(N( K_{l_x} ) ,\Z)$ and as 
$\bar{l}_1=d.\bar{l}_0$,  the  Mayer-Vietoris sequence  has the following shape: 
$$ ... \stackrel {\delta_2}{\longrightarrow}  H_1(  b( N( K_{l_x}) ),\Z)  \stackrel {\Delta_1}{\longrightarrow} \Z.  \bar{m} \oplus  \Z.  \bar{l}_0 \stackrel {i_1}{\longrightarrow} H_1 (L_X,\Z) \to 0$$
 where   $\Delta_1(\bar{l}_1)=( 0,d.\bar{l}_0)$ and  $\Delta_1(\bar{m})=(\bar{m},0)$. So,
  $H_1(L_X,\Z )$ is isomorphic to $\Z/d.\Z.$

\end{example}

We generalize the notion of $d$-curlings as follows.  Let  $k(D)$  be the $k$-pinched disc,   where $k>1$, quotient by identification  of their centrum of $k$ oriented  and ordered discs $D_i, 1\leq  i\leq k$. Let $c=c_1 \circ c_2  \circ  \dots  \circ c_n$  be a permutation of  the indices $\{1, \dots ,k\}$  given as the composition of $n$ disjoint cycles $c_j, 1\leq j \leq n,$  where $c_j$ is the  cycle $(1(j),...,i(j), ..., d_j(j))$ of order $d_j, 1\leq d_j$.  Let $\tilde{h}_c$ be an orientation preserving homeomorphism of  the disjoint union of  the $k$ discs  $D_i, 1\leq  i\leq k$, such that $\tilde{h}_c (D_i)=D_{c(i)}$,  $ \tilde{h}_c(0_i)=0_{c(i)}$ and $\tilde{h}_c^{d_j}$ restricted to $\coprod_{1\leq i \leq d_j} D_{i(j)}$ is the identity.  Let $h_c$ be the orientation preserving homeomorphism of $k(D)$ induces by $\tilde{h}_c$. By construction we have  $h_c(\tilde{0})= \tilde{0}$.  Up to isotopy the homeomorphisms $\tilde{h}_c$ and $h_c$ depend only on the permutation $c.$ So, the homeomorphism class of the mapping torus of $h_c$ acting on $k(D)$ depends only on $c.$ So, we can  state  the following definition. 

\newpage

\begin{definition}\label{deperm}
   
     A  {\bf singular pinched  solid torus  associated to the permutation $c$} is a topological space orientation preserving  homeomorphic to the  mapping torus  $T(k(D),c)$ of $h_c$ acting on $k(D)$:
 $$ T(k(D),c)=  \lbrack 0,1 \rbrack  \times k(D) / (1,x) \sim ( 0,h(x))$$
  The {\bf core } of $T(k(D),c)$ is the  oriented circle    $ l_0=   \lbrack  0,1 \rbrack \times  \tilde{0}  / (1,\tilde{0}) \sim   (0,\tilde{0})$.  A  homeomorphism between two singular pinched  solid tori is orientation preserving if it preserves the orientation of  $k(D)\setminus \{ \tilde{0} \}$ and the orientation of  the core $l_0$. 
  \end{definition}
  
 \begin{definition}\label{desheet}
 
 Let $T(k(D),c)$ be the singular pinched torus associated to the permutation $c$.{\bf{ A sheet of $T(k(D),c)$ }} is the closure of a connected component of $(T(k(D),c)\setminus l_0).$
 
 \end{definition}

\begin{remark}\label{remperm}  By construction 
 $( T(k(D),c) \setminus l_0)$
 has n connected components, each of them  is homeomorphic to  a torus minus its core $(S\times D)\setminus (S\times \{0\})$. The closures,  in $T(k(D),c) $,  of these  connected components,  are  the  mapping tori  of   the  homeomorphisms which permute cyclically  the   $d_j$ irreducible components
  of the  pinched discs $d_j(D)$. So,  Point 6 of Remark \ref{remcurl} implies that  a  sheet    $T(d_j(D),c_j), 1\leq j \leq n, $ of  $T(k(D),c) ,$   is a    $d_j$-curling.
 
 As  $k>1$, a homeomorphism of $T(k(D),c) $ preserves the core $l_0$ which is the set of the topologically singular points of $T(k(D),c) $. So the homeomorphism classes  of the  sheets  $T(d_j(D),c_j), \ 1\leq j\leq n ,$    of $T(k(D),c)$,    depend only on the permutation $c.$ The identifications point by point of the $n$ circles which are the $n $ cores $l_j$ of $T(d_j(D),c_j), \ 1\leq j\leq n $,  determine a quotient morphism 
  $$ \delta : (\coprod _{1\leq j\leq n} T(d_j(D),c_j)) \to   T(k(D),c) .$$ 
 Conclusion:  Up to homeomorphism  $T(k(D),c) $ can be obtained  by the composition of   the  $d_j$-curling morphisms  $q_j :(S\times D)\to {\cal {C}}_{d_j} , 1\leq j \leq n$,  defined on the disjoint union of $n$ solid tori,   followed by the identification  $\delta $ of their cores.

\end{remark}

\section{The topology of the normalization}

Let $(X,0) $ be a reduced surface germ, let $(\Sigma, 0)$ be its singular locus  and let 
 $\nu  :  (X',p') \rightarrow (X,p)$  be its  normalization. As in  Definition \ref{de1},  if $\sigma$ is  an irreducible component of $\Sigma $, let $\sigma '_j , 1\leq j \leq n(\sigma), $ be the  $n(\sigma)$  irreducible components of $\nu ^{-1} (\sigma)$,    and   let $d_j$ be the degree of $\nu $ restricted to $\sigma '_j $.  Moreover let
  $ k(\sigma)=:d_1+\dots +d_j+\dots +d_{n(\sigma )}$ be the total degree of $\nu$ above $\sigma$.\\
 Let $\Sigma_+$ be the union of the irreducible components $\sigma $ of  $\Sigma $ such that $k(\sigma )>1$. In $L_X$,  let $K_{\Sigma_+}$ be the link of $\Sigma_+$.  We choose  a compact regular neighbourhood $N(K_{\Sigma_+})$ of $K_{\Sigma_+}$. Let $E(K_{\Sigma_+})$ be the closure of $L_X \setminus N(K_{\Sigma_+})$. By definition  $E(K_{\Sigma_+})$ is the (compact) exterior of  $K_{\Sigma_+}.$

 \begin{lemma}\label{lem1}

 1.  The restriction of $\nu $ to $\nu ^{-1} (E(K_{\Sigma_+})) $ is an homeomorphism and  $(L_X\setminus K_{\Sigma_+})$ is a topological manifold.\\
 2. The link   $K_{\Sigma_+}$ is the set of the topologically singular points of $L_X$.\\
   3.  The homeomorphism class of $L_X$  determines  the homeomorphism class of
   $N(K_{\Sigma_+})$ and  $E(K_{\Sigma_+})$. \\
    4.  The number of connected components of $E(K_{\Sigma_+})$ is equal to the number of irreducible components of $(X,0).$

\end{lemma} 

{\it Proof:} \\ If $(X,0)$ has an isolated singular point at the origin, the normalization is bijective and as  the links  $L_X$  and $L_{X'}$ are  compact, $\nu$ restricted to $L_{X'}$ is a homeomorphism.

If $\Sigma $ is one dimensional,   let $K_{\Sigma}$ be the link of $\Sigma$.  In $L_X$ we choose  a compact regular neighbourhood $N(K_{\Sigma})$ of $K_{\Sigma }$. Let $E(K_{\Sigma })$ be the closure of $L_X \setminus N(K_{\Sigma })$. By definition  $E(K_{\Sigma })$ is the (compact) exterior of  $K_{\Sigma }.$ As $\nu $ restricted to $(X\setminus \Sigma)$ is bijective, $\nu $ restricted to the compact $\nu ^{-1} (E(K_{\Sigma})) $ is a homeomorphism. 

Let $\sigma $ be an irreducible component of $\Sigma $ and let $N(K_{\sigma})$ be the connected component of $N(K_{\Sigma })$ which contains the link $K_{\sigma}$. \\
 When $k(\sigma)=1$, $\nu $ restricted to  $\nu ^{-1}(N(K_{\sigma}))$ is a bijection. So,   the restriction of $\nu $ to $\nu ^{-1} (E(K_{\Sigma_+}) )$ is an homeomorphism.  Moreover $\nu$ restricted to $ \nu^{-1}( X \setminus \Sigma )$ is an analytic isomorphism. So, $(L_X\setminus K_{\Sigma_+})$ is a topological manifold. This ends the proof of Statement 1.

  If $k(\sigma)>1$,  $\nu $ restricted  to $\nu ^{-1} (K_{\sigma}) $ is not injective. Let $p$ be a point of $K_{\sigma}$.  The number of the irreducible components  $\sigma'_j$ of   $\nu ^{-1} (\sigma)$ is denoted  $n(\sigma )$. So,      $\nu ^{-1} (\sigma )= \cup _{1\leq j \leq n(\sigma)}\sigma '_j  $.    Let $d_j$ be  the degree of $\nu $ restricted to $\sigma'_j$. The intersection   $\nu ^{-1}(p)\cap \sigma '_j$  has $d_j$ points $ \{ p_{i(j)},  1\leq i \leq d_j  \}$. As $(X',p')$ is normal,   $(X'\setminus p')$ is smooth  and $(\sigma'_j \setminus p')$ is a smooth curve at   any point $z_j \in (\sigma'_j \setminus p')$.  
 In $(X'\setminus p')$,  we can choose at the  points $p_{i(j)}$, a smooth germ  of curve  $(\gamma _{i(j)},p_{i(j)})$ which cuts $\sigma'_j$ transversally at $p_{i(j)}$ and such that $ D'_{i(j)}= \nu^{-1}(N(K_{\sigma})) \cap \gamma _{i(j)}$  is a disc centered at $p_{i(j)}$. Let $D_{i(j)}$ be $\nu(D'_{i(j)})$. By construction $p$ is the common center of the topological discs $D_{i(j)} $. So,
   $ (\cup  _{1\leq j \leq  n(\sigma )} (\cup_ {1  \leq i \leq d_j} D_{i(j)}) )$ 
    is a $k(\sigma )$-pinched disc centered at $p$. As $k(\sigma)>1$, $L_X$ is not a topological  manifold at $p$. This ends the proof of Statement   2.

Statements  1 and  2 imply  that  $(K_{\Sigma_+})$ is the set of the topologically  singular points of $L_X$. It implies  3.

  The number $r$ of irreducible components of  $(X,0)$ is  equal to the number of connected components of $L_{X'}$. But, $L_{X'}$ and $\nu ^{-1} (E(K_{\Sigma_+}) )$ have the same number of connected components since $( L_{X'} \setminus \nu ^{-1} (E(K_{\Sigma_+}) ))$ is a regular neighbourhood of the differential link $\nu ^{-1} (K_{\Sigma_+}) $. Statement 1 implies that  $r$ is also the number of connected components  of  $E( K_{\Sigma_+})$.  This proves 4.\\
  {\it End of proof.}

\begin{lemma}\label{lem2} 

 Let  $\sigma$ be an irreducible component of $\Sigma_+ $,   let $K_{\sigma}$ be the link of $\sigma $ in $L_X$. We choose,  in $L_X$,  a compact regular neighbourhood $N(K_{\sigma})$ of $K_{\sigma}$.  
The link    $K_{\sigma}$ is a deformation retract of  $N(K_{\sigma})$.  If  $l_{\sigma } $ is  the homotopy class of  $K_{\sigma}$  in $\pi _{1} (N(K_{\sigma }) )$,   then   $\pi _1 (N(K_{\sigma}) = \Z  .  l_{\sigma}$.
We have:

1.  The tubular neighbourhood   $ \nu ^{-1}(N(K_{\sigma}))$ of $ \nu ^{-1}(K_{\sigma})$ is the disjoint union of $n(\sigma )$ solid tori\\ $T'_j,  1\leq j \leq n(\sigma )$,  and 
the boundary of  $N(K_{\sigma})$ is the disjoint union of $n(\sigma )$  tori.

2.  Let $c$ be the permutation of $k(\sigma )$ elements which is the composition  of $n(\sigma)$ disjoint cycles $c_j$ of order $d_j$.  Then  $N(K_{\sigma})$ is homeomorphic to a singular pinched torus $T(k(\sigma)(D),c)$ wich   has $n(\sigma)$ sheets $T_j$ where  
  $T_j=\nu(T'_j) , 1\leq j \leq n(\sigma)$.

3.   On  each connected component $\tau_j$ of the boundary of $N(K_{\sigma})$, the homeomorphism class of $N(K_{\sigma})$ determines a  unique (up to isotopy)  meridian curves $m_j$.  If $l_j$ is a  parallel on  $\tau _j$ the homotopy class of $l_j$ in $\pi _{1} (N(K_{\sigma }) )$ is equal to $d_j.l_{\sigma }$.

 \end{lemma}

{\it Proof of Lemma \ref{lem2}:} \\
 The link $L_{X'}$ of the normalization  $\nu  :  (X',p') \rightarrow (X,p)$ of $(X,0)$ is a three dimensional Waldhausen graph manifold. Let $\sigma $ is an irreducible component of $\Sigma_+ $. 
  In $L_{X'}$,  $\nu^{-1}(K_{\sigma})$ is a differentiable one dimensional link with $n(\sigma)$ irreducible components  $K_{\sigma'_j}, 1\leq j \leq n(\sigma )$.   But ,  $\nu^{-1}(N(K_{\sigma}))$  is a regular compact neighbourhood  of the link $( \coprod _{ 1\leq j \leq n(\sigma )} K_{\sigma'_j})$. So,  th   $\nu^{-1}(N(K_{\sigma}))$ is the disjoint union of $n(\sigma )$ solid  tori that we denote by  $T'_j , 1\leq j \leq n(\sigma )$.  Moreover   let $\tau'_j$ be the boundary of $T'_j$.  As $\nu $ restricted to the boundary of  $\nu^{-1}(N(K_{\sigma}))$ is a homeomorphism, the boundary of  $N(K_{\sigma})$ is the disjoint union of  the $n(\sigma )$ tori $\tau_j=\nu (\tau'_j)$.  Statement 1 is proved.

   As in the proof of Lemma \ref{lem1}, we consider a point $p\in K_{\sigma}$ and  $d_j$ meridian discs $D'_{i(j)}$  of $T'_j$  centered at the $d_j$  points of $\nu^{-1}(p) \cap \sigma'_j $.  We equip $T'_j$ with a trivial fibration in  oriented circles  which has $K_{\sigma'_j}$ as central  fiber. The first return homeomorphism $h'_j: (\coprod_i D'_{i(j)}) \to (\coprod_i D'_{i(j)}) $,  along the  chosen circles  permutes cyclically the $d_j$ discs $D'_{i(j)}$ and $(h'_j)^{d_j}$ is the identity. So,  $h_j$ provides  a $d_j$-cycle $c_j$. 
   
    The direct image by $\nu$,  of the fibration  of $T'_j $ in circles,  equip $T_j=\nu (T'_j)$ with a foliation in oriented circles which has $K_{\sigma}$  as  singular  leave. But $ \nu (\coprod_i D'_{i(j)})$ is a $d_j$-pinched disc with origin $p\in K_{\sigma}$ and  irreducible components $D_{i(j)} =\nu (D'_{i(j)})$. Let $h_j$ be the homeomorphism defined on the $d_j$-pinched disc  $(\cup_i D_{i(j)})$ by the first return along the given circles. By construction, $h_j(p)=p$ and $h_j$ permutes cyclically the $d_j$ irreducible components of  the $d_j$-pinched disc   $(\cup_i D_{i(j)})$.  So, $T_j$ is the mapping torus of  $h_j$ acting on  the $d_j$-pinched disc $(\cup_i D_{i(j)})$.  By  Point 6 of Remark \ref{remcurl},  $T_j$ is a $d_j$-curling. 
    
    But $h_j(p)=p$ for all $j, 1\leq j\leq n(\sigma)$. So, if  we consider $h(z)=h_j(z)$ for all  $z\in (\cup_i D_{i(j)}))$, $h$ is  a well defined  homeomorphism 
    $$h:(\cup_{1\leq j\leq n(\sigma)} (\cup_i D_{i(j)})) \to  (\cup_{1\leq j\leq n(\sigma)} (\cup_i D_{i(j)})) .$$ 
     By construction, $h$ induces on  the $k(\sigma)$-pinched disc $(\cup_{1\leq j\leq n(\sigma)} (\cup_i D_{i(j)}))$ a permutation $c$,  of its irreducible components, which is the  composition of the disjoint  cycles $c_j$. Then,   $N(K_{\sigma})$ is homeomorphic to the singular pinched torus $T(k(\sigma)(D),c)$. But 
      $T_j, 1\leq j\leq n(\sigma)$,  is the closure of a connected components of 
      $(N(K_{\sigma}) \setminus K_{\sigma})$. By definition the family of the    $d_j$-curlings  $ \{ {\cal {C}}_{d_j}=T_j, 1\leq j\leq n(\sigma) \},$ is  the  family of the  sheets of $N(K_{\sigma})$.

 A  pinched disc can be retracted by deformation onto its center. Such a retraction by  deformation can be extended in a retraction by deformation of  $T(k(\sigma)(D),c)$ onto  its core. By  2,     $N(K_{\sigma})$ is homeomorphic to a singular pinched torus $T(k(\sigma)(D),c)$ and $K_{\sigma}$ is its core. Then, we can retract $N(K_{\sigma})$ by deformation onto its core $K_{\sigma}$.
   So,  $\pi _1 (N(K_{\sigma}) )= \Z  .  l_{\sigma}$ where  $l_{\sigma } $  is  the homotopy class of  $K_{\sigma}$  in $\pi _{1} (N(K_{\sigma }) )$.  
 Moreover each connected component $\tau_j$ of  the boundary of $N(K_{\sigma})$ is the boundary of     the $d_j$-curling  $T_j$ which is a sheet of  $N(K_{\sigma})$.  On $\tau_j$ we choose a meridian curve  $m_j$  and a parallel  curve  $l_j$  as defined in \ref {decurl}.   As $T_j$ is a $d_j$-curling,  $m_j$ is unique up to isotopy by   1 of Remark \ref{remcurl}. By  2 of Remark \ref{remcurl}, the class of   $l_j$ in $\pi _1 (N(K_{\sigma})$ is equal to $d_j.l_{\sigma}$. \\
 {\it  End of proof of Lemma \ref{lem2}.}

Now we are ready to prove the main theorem of this paper.

\begin{theorem}\label{th1}
 Let $(X,0)$ be a reduced surface germ. The homeomorphism class of $L_X$  determines the homeomorphism class of the link $L_{X'}$ of the normalization of $(X,0)$.
\end{theorem} 

{\it  Proof:} \\ If $L_X$ is a topological manifold, Statements 1 and 2 of  Lemma \ref{lem1}  state that $K_{\Sigma_+}$ is empty. So,  $E(K_{\Sigma_+}) = L_X$ and $\nu $  is a homeomorphism. When $L_X$ is a topological manifold  the theorem is trivial.

If $L_X$ is not a topological manifold, let  ${\cal {K}}$ be  the set of its singular points. By Statement  2 of  Lemma \ref{lem1}, we know that ${\cal {K}}=K_{\Sigma_+}$  is a disjoint union of circles.
  Let $N{\cal {K}} )$ be a regular compact neighbourhood of $ {\cal {K}}$. 
 By definition,  the exterior $E({\cal {K}})$ of ${\cal {K}}$ is the closure of  $( L_X \setminus N({\cal {K}} ))$.

Let $K$ be a connected component  of ${\cal {K}}$ and let $N(K)$ be the connected component of $N({\cal {K}})$ which contains  $K$. There exist a irreducible component $\sigma $ of $\Sigma_+$ such that  $K=K_{\sigma }$. By Lemma \ref{lem2}, $N(K)$ is homeomorphic to a singular pinched torus $T(k(\sigma)(D),c)$ wich  has $n(\sigma)$ sheets  homeomorphic to $d_j$-curlings  ${\cal {C}}_{d_j} , 1\leq j \leq n(\sigma)$. As explained Section 2, the integers $d_j$ and $n(\sigma)$, the permutation $c$ and the homeomorphism class of  the family of the $d_j$-curlings  ${\cal {C}}_{d_j} $,  sheets of $N(K)$,  depend only on the homeomorphic class of $N(K)$.
In particular the boundary $b(N(K))$ has   $n(\sigma )$  tori $\tau_j,  1\leq j\leq n(\sigma),$ as connected components.\\
 As defined in \ref{demer} and shown in Statement  1 of Remark \ref{remcurl}, on     each torus  $\tau_j$ we have a well defined meridian curve $m_j$ associated to the corresponding sheet of $N(K)$. It is the key point of this proof. The existence of  well defined meridian curves  of  $N({\cal {K}})$ on  each torus of the boundary of the exterior $E({\cal {K}})$ of ${\cal {K}}$ allows us to perform Dehn fillings. As justify  below, to take  $E({\cal {K}})$  and to close   it  by performing    Dehn fillings associated to  the given meridian curves  of  $N({\cal {K}})$,  produce a closed manifold  homeomorphic to $L_{X'}$. Let us be more  precise.

{\bf { The  Dehn filling construction:}} \\   
Let $T$ be a solid torus given with a meridian disc $D$ and let  $m_T $ be the   boundary of $D$. By definition  $m_T$ is a meridian curve on the boundary of $T$.  Let $U(D)$ be  a compact  regular neighbourhood of $D$ in $T$ and let  $B$ be the closure of $T\setminus U(D)$.  By construction $B$  is a 3-dimensional ball.  In the   boundary   $b(T)$ of $T$, the closure of   the complement of the annulus $ U(m_T)=U(D)\cap b(T)$ is also an annulus $E(m_T )\subset b(B)$.\\
  On  the other hand, we  suppose that a torus  $\tau$ is  a   boundary component of an oriented compact three-dimensional  manifold  $M.$ 
 Let  $\gamma  $ be an oriented essential simple closed curve   on  $\tau $. So,  $\tau $ is the union of two annuli, $U(\gamma )$,  a compact regular neighbourhood  of $\gamma$,   and the closure $E(\gamma)$ of $\tau \setminus U(\gamma)$. There is a unique way to glue  $T$  to $M$  by an orientation reversing  homeomorphism   between  the boundary of  $T$  and  $\tau $ which send $m_T$ to $\gamma .$  
  Indeed,  the  gluing  of $U(m_T)$ onto $U(\gamma)$ determines a union $M'$ between  $U(D)$ and $ M$.  This gluing extends to  the gluing of $E(m_T)$ onto $E(\gamma)$ which determines a union between $B$ and $M'$. \\ So, the result of such a gluing  is  unique up to orientation preserving homeomorphism and it is called  the {\bf  Dehn filling } of $M$ associated to $\gamma $. \\
  One can find a presentation of the Dehn filling construction in S. Boyer \cite{boyer}. 
 
 The topology of the   link $L_X$ determines   the exterior  $E({\cal {K}})$ of the  singular locus ${\cal {K}}$ of $L_X$  and also the  well defined meridian curves  of $N({\cal {K}})$ on each connected component of the boundary of $E({\cal {K}})$.  Let $\sigma_i, 1 \leq i \leq r$ be the $r$ irreducible components of $\Sigma_+$.  So, the boundary of   $E({\cal {K}})$ has $n=\sum _{1\leq i \leq r} \  n(\sigma _i)$  connected components.  Let $T$ be a solid torus and let $m_T$  be a meridian curve of $T$. Let $\tau $ be  one connected component of the boundary of $E({\cal {K}})$ given with its  already chosen   curve $m$ which is a meridian curve of $N({\cal {K}}) $. By Remark \ref{remcurl},  $m$ is an essential  simple closed curve on  $\tau$.  We glue $T$ to   $E({\cal {K}})$ with the help of  an orientation reversing homeomorphism 
   $$f:  b(T)\to \tau $$ 
  defined on the boundary $b(T)$ of $T$  such that  $f(m_T)=m$.

  We perform such a  Dehn filling associated to $m$  on each of  the $n$ connected components of the boundary of $E({\cal {K}})$. So, we obtain a closed 3-dimensional Waldhausen graph manifold $L$.
 
 But,  $\nu $ restricted to $E'=\nu^{-1}( E({\cal {K}}))$  is a homeomorphism. Moreover, 
 $\nu^{-1} (N({\cal {K}}))$ is a tubular neighbourhood of the  differential link $\nu^{-1} ({\cal {K}})$ which has  $n=\sum_{1\leq i \leq r} \  n(\sigma _i)$ connected components.  So, $\nu^{-1} (N({\cal {K}}))$ is a disjoint union of $n$ solid tori. As in the proof of Lemma \ref{lem2}, let   $T'_j$ be one of these solid tori. Then,  $T_j=\nu (T'_j)$ is a sheet of $N(K_{\sigma})$ where $\sigma $ is an irreducible component of  $\Sigma _+.$ But $T_j$ is a $d_j$-curling and $\nu $ restricted to $T'_j$ is a quotient morphism associated to this $d_j$-curling.  Let $m_j$ be the chosen  meridian on $T_j$. Definition \ref{decurl} implies that $\nu ^{-1}(m_j)$ is a meridian curve of $T'_j$. By the  unicity of the Dehn filling construction, there exists an orientation preserving  homeomorphism between $L_{X'} $ and $L$.

Conclusion: The construction of    $L$   only  depends  on the topology of  $L_X$ and $L$   is homeomorphic to  the link $L_{X'}$ of the normalization  $(X',p')$ of $(X,0)$.  So, Theorem \ref{th1} is proved.

{\it End of proof.}

\newpage

\section{\bf  Surface germs with simply connected links}

This section is devoted to the proof of the following proposition.

   \begin{proposition}\label{pr1} Let $(X,0)$ be an irreducible surface germ. If the link $L_X$ of $(X,0)$ is simply connected then the normalization $\nu : (X',p') \to (X,0)$ is a homeomorphism and $(X',p')$ is smooth at $p'$. In particular,  the normalization is the good minimal resolution of  $(X,0)$.
 
 \end{proposition}

{\it  Proof:}  \\  By  Lemma \ref{lem1} (or Proposition 3.12 in \cite{mi}),  if  $L_X$ is a topological manifold   the normalization $\nu  :  (X',p') \rightarrow (X,p)$ is a homeomorphism.  Then the link $L_{X'}$ is also simply connected and by Mumford's theorem \cite{mu} $(X',p')$ is smooth at $p'$.
 
 Now, we  suppose that $L_X$ is not a topological manifold.  Then,  the following two statements I and II prove that $L_X$ is not simply connected.\\ As before,    $\Sigma_+$ is the union of the irreducible components  $\sigma$  of the singular locus of $(X,0)$, which have  a total degree,  $ k(\sigma)= d_1+\dots +d_j+\dots +d_{n(\sigma )}$,  stricly greater than one.
  By Lemma \ref{lem1},  if $L_X$ is not a topological manifold,   $\Sigma_+$ has at least one irreducible component $\sigma .$

{\bf { Statement I }} If there exists an irreducible component  $\sigma $  of $\Sigma _+$ with $n(\sigma) >1$, then  the rank $r$ of $H_1(L_X, \Z)$  is greater than  or equal to $  (n(\sigma )-1)$, in particular $H_1(L_X, \Z)$ has infinite order.

{\it Proof of Statement I.} Let $ \Sigma_+ = (\sigma \cup_{1\leq i \leq r } \sigma_i ),$ be the decomposition of $\Sigma_+$ as the  union of its  irreducible components.  As $(X,0)$ is irreducible,  $E( K_{\Sigma_+} ) $ is connected  by Lemma \ref{lem1}. Then, $ E( K_{\sigma})=E( K_{\Sigma_+}) \cup _{1\leq i \leq r} N(K_{\sigma _i} ) $  and $L_X$  are also  connected. But $N(K_{\sigma })$ which is a singular pinched torus (Lemma \ref{lem2}) is connected with $n(\sigma)>1$ boundary components.  We consider the Mayer-Vietoris  exact sequence associated to the decomposition of $L_X$ as  the union $  E( K_{\sigma})  \cup   N( K_{\sigma}) .$
$$ ...\to  H_1(L_X, \Z) \stackrel {\delta_1}{\longrightarrow}  H_0(  E( K_{\sigma}) \cap  N( K_{\sigma}) ,\Z)  \stackrel {\Delta_0}{\longrightarrow}  H_0( E( K_{\sigma}), \Z) \oplus  H_0(  N( K_{\sigma}) ,\Z)  \stackrel {i_0}{\longrightarrow}  H_0 (L_X,\Z) \to 0$$

 But  $E( K_{\sigma}) \cap  N( K_{\sigma})$  is the disjoint union of $n(\sigma)$ disjoint tori. So, the rank of $ H_0(  E( K_{\sigma}) \cap  N( K_{\sigma}) ,\Z)$ is equal to $n(\sigma).$ Since $\sigma $ is irreducible $ H_0(  N( K_{\sigma}) ,\Z)$  has rank one.  Since $(X,0)$ is irreducible,   $ H_0( E( K_{\sigma}),\Z) $  and  $ H_0 (L_X,\Z) $ have rank one. So,  the rank of $ Ker (\Delta_0)=\delta_1 (H_1(L_X, \Z))$ is equal to $(n(\sigma)-1)$. This ends the proof of Statement I.

{\bf { Statement II }}  If there exists an irreducible component  $\sigma $  of $\Sigma _+$ with $n(\sigma) =1 $ and $k(\sigma)=d>1$ the order of $H_1(L_X, \Z)$ is at least  $d$.

 {\it Proof of Statement II.} By Lemma \ref{lem2}, if $n(\sigma)=1$ and $d>1$,  $ N( K_{\sigma}) $ is a $d$-curling and the boundary of $ N( K_{\sigma}) $ is a torus $\tau$. By Remark \ref{remcurl}, $\tau$ is  given with a meridian curve $m$ and a parallel curve $l$.  Let $\bar{m}$, $\bar{l}$ and $l_{\sigma}$ be the classes of $m$,  $l$ and $K_{\sigma}$,  in $H_1(  N( K_{\sigma}) ,\Z) $.  Moreover, we have   $ H_1(  N( K_{\sigma}) ,\Z)   =\Z .l_{\sigma}$ and  $ \bar{l}=d.l_{\sigma}$. 
We consider the Mayer-Vietoris  exact sequence associated to the decomposition of $L_X$  as  the union $  E( K_{\sigma})  \cup   N( K_{\sigma}) .$
$$ ...\to  H_2(L_X, \Z) \stackrel {\delta_2}{\longrightarrow}  H_1(  E( K_{\sigma} ) \cap  N( K_{\sigma}) ,\Z)  \stackrel {\Delta_1}{\longrightarrow}  H_1( E( K_{\sigma}) ,\Z ) \oplus  H_1(  N( K_{\sigma}) ,\Z)) \stackrel {i_1}{\longrightarrow} H_1 (L_X,\Z) \to ...$$
As  $E( K_{\sigma}) \cap  N( K_{\sigma})=\tau$, the image of  $\Delta_1 $ is generated by $\Delta_1(\bar{m})=(x,0)$ and $\Delta_1(\bar{l})=(y,d.l_{\sigma})$ where $x$ and $y$ are in $H_1( E( K_{\sigma}),\Z )$. So,  the image of $\Delta _1$ is included in  $ H_1( E( K_{\sigma}),\Z ) \oplus  \Z \ d.l_{\sigma}.$ It implies that  the order of the cokernel  of $\Delta_1$ is at least  $d$. This ends the proof of Statement II.  \\
 The  two statements above imply  Proposition \ref{pr1}.\\ {\it End of proof.}

\vskip.1in

\end{document}